\documentclass{amsart}

\usepackage{amssymb}

\usepackage{amsmath,  amsopn, amsfonts}

\newtheorem{theorem}{Theorem}[section]
\newtheorem{lemma}[theorem]{Lemma}
\newtheorem{corollary}[theorem]{Corollary}

\theoremstyle{definition}
\newtheorem{definition}[theorem]{Definition}

\newtheorem{proposition}[theorem]{Proposition}

\theoremstyle{remark}
\newtheorem{remark}[theorem]{Remark}

\numberwithin{equation}{section}

\renewcommand\Re{\mathrm{Re}}
\renewcommand\Im{\mathrm{Im}}
\newcommand\R{{\mathbb R}}
\newcommand\N{{\mathbb N}}
\newcommand\Z{{\mathbb Z}}

\newcommand\eqdef {\buildrel\hbox{{\footnotesize def}}\over =}

\begin{document}

\title{Square function and heat flow estimates on domains}

\author{O.Ivanovici}
\email{oana.ivanovici@math.u-psud.fr}
\address{Universit{\'e} Paris Sud,
Math{\'e}matiques,
B{\^a}t 430, 91405 Orsay Cedex}

\author{F.Planchon}
 \email{fab@math.univ-paris13.fr}
\address{Universit{\'e} Paris 13,
L.A.G.A.,UMR 7539, Institut Galil\'ee,
99 avenue J.B,Cl\'ement, F-93430 Villetaneuse}
\thanks{}

\subjclass[2000]{Primary 35J25,58G11}

\date{}

\begin{abstract}
The first purpose of this note is to provide a proof of the usual
square function estimate on $L^p(\Omega)$. It turns out to follow
directly from a generic Mikhlin multiplier theorem obtained by Alexopoulos,
which mostly relies on Gaussian bounds on the heat kernel. We also
provide a simple proof of a weaker
version of the square function estimate, which is enough in most
instances involving dispersive PDEs. Moreover, we obtain, by a
relatively simple integration by parts, several useful $L^p(\Omega;H)$
bounds for the derivatives of the heat flow with values in a given
Hilbert space $H$.
\end{abstract}

\maketitle

\bibliographystyle{plain}

\section{Introduction}
Let  $\Omega$ be a domain in $\mathbb{R}^{n}$,  $n\geq 2$, with smooth boundary $\partial\Omega$. Let $\Delta_{D}$ denote the Laplace operator on $\Omega$ with Dirichlet boundary conditions, acting on $L^{2}(\Omega)$, with domain $H^{2}(\Omega)\cap H^{1}_{0}(\Omega)$.

The first result reads as follows:
\begin{theorem}\label{theoremtriebel}
Let $f\in C^{\infty}(\Omega)$ and $\Psi\in C^{\infty}_{0}(\mathbb{R}^{*})$ such that
\begin{equation}\label{unu}
\sum_{j\in\Z}\Psi(2^{-2j}\lambda)=1,\quad \lambda\in\mathbb{R}.
\end{equation}
Then for all $p\in (1,\infty)$ we have
\begin{equation}\label{ib}
\|f\|_{L^{p}(\Omega)}\approx C_{p}\Big\|\Big(\sum_{j\in\Z}|\Psi(-2^{-2j}\Delta_{D})f|^{2}\Big)^{1/2}\Big\|_{L^{p}(\Omega)},
\end{equation}
where the operator $\Psi(-2^{-2j}\Delta_{D})$ is defined by
\eqref{integralpsi} below.
\end{theorem}
Readers who are familiar with functional spaces'theory will have
recognized the equivalence $\dot F^{0,2}_p\approx L^p$, where the Triebel-Lizorkin space is defined using the right hand-side
of \eqref{ib}
as a norm. In other words, $L^p(\Omega)$ and the
Triebel-Lizorkin space $\dot F^{0,2}_p(\Omega)$ coincide. Such an
equivalence (and much more !) is proven in
\cite{Triebel1,Triebel2,Triebel3}, though one has to reconstruct it
from several different sections (functional spaces are defined
differently, only the inhomogeneous ones are treated, among other
things). As such, the casual user with mostly a PDE background might
find it difficult to reconstruct the argument for his own sake without
digesting the whole theory. It turns out that the proof of \eqref{ib}
follows directly from the classical argument (in $\R^n$) involving
Rademacher functions, provided that an appropriate Mikhlin-H\"ormander
multiplier theorem is available. We will provide details below.

A weaker version of Theorem \ref{theoremtriebel} is often used in the
context of dispersive PDEs:
\begin{theorem}\label{theorembesov}
Let $f\in C^{\infty}(\Omega)$, 
then for all $p\in [2,\infty)$ we have
\begin{equation}\label{ibp}
\|f\|_{L^{p}(\Omega)}\leq C_{p}|\Big(\sum_{j\in\Z}\|\Psi(-2^{-2j}\Delta_{D})f\|^2_{L^{p}(\Omega)}\Big)^{1/2}.
\end{equation}
\end{theorem}

 The second part of the present note aims at
giving a self-contained proof of \eqref{ibp} , with ``acceptable'' black boxes,
namely complex interpolation and spectral calculus. In fact, if one accepts to replace the spectral
localization by the heat flow, the proof can be made entirely
self-contained, relying only on integration by parts. Our strategy to
prove Theorem \ref{theorembesov} is indeed to reduce matters to an
estimate involving the heat flow, by proving almost orthogonality
between spectral projectors and heat flow localization; this only
requires basic parabolic estimates in $L^p(\Omega)$, together with a
little help from spectral calculus. 

\begin{remark}
For compact manifolds without boundaries, one may find a direct proof
of \eqref{ibp} (with $\Delta_D$ replaced by the Laplace-Beltrami
operator) in \cite{BGT1}, which proceeds by reduction to the $\R^n$
case using standard pseudo-differential calculus. Our elementary approach
provides an alternative direct proof. However, the true square
function bound \eqref{ib} holds on such manifolds, as one has a
Mikhlin-H\"ormander theorem from \cite{SeeSog}.
\end{remark}
\begin{remark}
  One can also adapt all proofs to the case of Neumann boundary
  conditions, provided special care is taken of the zero frequency
  (note that on an exterior domain, a decay condition at infinity
  solves the issue). The Gaussian bound which is required later holds in
  the Neumann case, see \cite{DavHeat,Daners}.
\end{remark}
\begin{remark}
  As mentioned before, Theorem \ref{theorembesov} is useful, among other things, when
  dealing with $L^p$ estimates for wave or dispersive evolution
  equations. For such equations, one naturally considers initial data
  in Sobolev spaces, and spectral localization conveniently reduces
  matters to data in $L^2$, and helps with finite speed of propagation
  arguments. One however wants to sum eventually over
  all frequencies in $l^2$, if possible without loss. Recent examples on
  domains may be
  found in \cite{OanaExtSchrod} or \cite{PV}, as well as in \cite{GillesEvian}.
\end{remark}
We now state estimates involving directly the heat flow, which will be
proved by direct arguments. It should be noted that for nonlinear
applications, it is quite convenient to have bounds on derivatives of
spectral multipliers, and such bounds do not follow immediately from
the multiplier theorem from \cite{Alex1}. We consider the linear heat
equation on $\Omega$ with Dirichlet boundary conditions and initial
data $f$
\begin{equation}\label{heateq}
\partial_{t}u-\Delta_{D}u=0, \ \text{on} \ \Omega\times\mathbb{R}_{+};\quad
u|_{t=0}=f\in C^{\infty}(\Omega);\quad u|_{\partial\Omega}=0.
\end{equation}
We denote the solution $u(t,x)=S(t)f(x)$, where we set
$S(t)=e^{t\Delta_{D}}$.  For the sake of simplicity $\Delta_{D}$ has
constant coefficients, but the same method applies in the case when
the coefficients belong to a bounded set of $C^{\infty}$ and the
principal part is uniformly elliptic (one may lower the regularity
requirements on both the coefficients and the boundary, and a nice
feature of the proofs which follow is that counting derivatives is
relatively straightforward).

Let us define two operators which are suitable heat flow versions of
$\Psi(-2^{-2j}\Delta_D)$:
\begin{equation}
  \label{eq:defQ}
  Q_t =\sqrt t \nabla S(t) \,\,\text{ and }\,\, \mathrm{\bf
    Q}_{t}\eqdef t\partial_t S(t).
\end{equation}
\begin{theorem}\label{propgrad}
Let $1< p<+\infty$, then we have
\begin{equation}
  \label{eq:quadheat}
  \|f\|_{L^{p}(\Omega)}\approx c_{,p}\|\Big(\int_{0}^{\infty}|Q_t
  f|^{2} \frac{dt} t \Big)^{1/2}\|_{L^{p}(\Omega)},
\end{equation}
which implies, for $p\in [2,+\infty)$,
\begin{equation}
  \label{eq:quadbheat}
  \|f\|_{L^{p}(\Omega)}  \leq C_{p}\Big(\int_{0}^{\infty}\|Q_t
  f\|^{2}_{L^{p}(\Omega)}\frac{dt} t \Big)^{1/2},
\end{equation}
and $Q_t$ may be replaced by $\mathrm{\bf Q}_t$ in both statements.
\end{theorem}
Notice that there is no difficulty to define $Q_t f$ or $\mathrm{\bf
  Q}_t f$ as distributional derivatives for $f\in L^p(\Omega)$, while
simply defining $\Psi(-2^{-2j}\Delta_D)$ on $L^p(\Omega)$ is already
a non trivial task. The purpose of the next Proposition is to prove that
both operators are in fact bounded on $L^p(\Omega)$.
\begin{proposition}\label{lemestvect}
Let $1<p<+\infty$. The operators $Q_{t}$, $\mathrm{\bf Q}_t$ are
bounded on $L^{p}(\Omega)$, uniformly in $t\geq 0$. Moreover
$\mathrm{\bf Q}_t$ is bounded on $L^1(\Omega)$ and $L^\infty(\Omega)$.
\end{proposition}
For practical applications, one may need a vector valued version of
Proposition \ref{lemest}. Let us consider now $u=(u_{l})_{l\in\{1,..,N\}}$ for $N\geq
2$, where each $u_{l}$ solves \eqref{heateq} with Dirichlet condition
and initial data $f_{l}$. Let $H$ be the Hilbert space with norm
$\|u\|^2_H=\sum_l |u_l|^2$, and $L^p(\Omega;H)$ the Hilbert valued
Lebesgue space. Then we have
\begin{proposition}\label{lemest}
Let $1<p<+\infty$. The operators $Q_{t}$, $\mathrm{\bf Q}_t$ are
bounded on $L^{p}(\Omega;H)$, uniformly in $t\geq 0$ and $N$. Moreover
$\mathrm{\bf Q}_t$ is bounded on $L^1(\Omega;H)$ and $L^\infty(\Omega;H)$.
\end{proposition}
\begin{remark}
One may therefore extend the finite dimensional case to any separable
Hilbert space. The typical setting would be to consider the solution $u$ to the
  heat equation with initial data $f(x,\theta)\in
  L^2_\theta=H$. Notice that the Hilbert valued bound does not follow
  from the previous scalar bound; however the argument is essentially
  the same, replacing $|\cdot|$ norms by Hilbert norms.
\end{remark}
\begin{remark}
  A straightforward consequence of Propositions \ref{lemest} and \ref{lemestvect} is
  that the Riesz transforms $\partial_j (-\Delta_D)^{-\frac 1 2}$ are
  continuous on Besov spaces defined by the RHS of \eqref{eq:quadbheat};
  these spaces are equivalent to the ones defined by the RHS of
  \eqref{ibp}, see Remark \ref{eqBesov} later on.
\end{remark}
Alternatively, one can derive all the (scalar, at least) results on the heat flow from
adapting to the domain case the theory which ultimately led to the
proof of the Kato conjecture (\cite{AT,AMT}). Such a possible
development is pointed out by P. Auscher in \cite{Auscher} (chap. 7,
p. 66) and was
originally our starting point; eventually we were led to the
elementary approach we present here, but we provide a sketch of
an alternate proof in the next  remark, which was kindly
outlined to us by Pascal Auscher. 
\begin{remark}
\label{rempascal}
The main drawback from \eqref{eq:quadbheat} is the presence of $\nabla
S(t)$ on the right hand-side: one is leaving the functional calculus of
$\Delta_D$, and in fact for domains with Lipschitz boundaries the operator $\nabla
S(t)$ may not even be bounded. As such, a suitable alternative is to
replace $\nabla S(t)$ by $\sqrt{ \partial_t} S(t)$. Then the
square function estimate may be obtained following \cite{Auscher} as
follows:
\begin{itemize}
\item prove that the associated square function in time is bounded by
  the $L^p$ norm, for all $1<p\leq 2$, essentially following step 3 in
  chapter 6, page 55 in
  \cite{Auscher}. This requires very little on the semi-group, and
  Gaussian bounds on $S(t)$ and $\partial_t S(t)$ (\cite{DaviesJOT})
  are more than enough to apply the weak $(1,1)$ criterion from
  \cite{Auscher} (Theorem 1.1, chapter 1). Moreover, the argument can be extended to domains
  with Lipschitz boundaries, assuming the Laplacian is defined
  through the associated Dirichlet form;
\item by duality, we get the square function bound for $p>2$ (step 5, page 56 in
  \cite{Auscher});
\item from now on one proceeds as in the remaining part of our paper
  to obtain the bound with spectral localization, and almost
  orthogonality \eqref{estkj} is even easier because we stay in the
  functional calculus. One has, however, to be careful if one is
  willing to extend this last step to Lipschitz boundaries, as this
  would most likely require additional estimates on the resolvent
  to deal with the $\Delta_j$.
\end{itemize}
\end{remark}
\section{From a Mikhlin multiplier theorem to the square function}
The following ``Fourier multiplier'' theorem is obtained in
\cite{Alex1} under very weak hypothesis on the underlying manifold
(see also \cite{Alex2} for a specific application to Markov chains,
and \cite{DOS} for a version closer to the sharp H\"ormander's multiplier
theorem, under suitable additional hypothesis, all of which are
verified on domains). For $m\in L^\infty(\R^+)$, one usually defines the
operator $m(-\Delta_D)$ on $L^2(\Omega)$ through the spectral measure
$dE_\lambda$:
\begin{equation}
  \label{eq:sM}
  m(-\Delta_D)=\int_0^{+\infty} \,dE_\lambda,
\end{equation}
and $m(-\Delta_D)$ is bounded on $L^2$.
\begin{remark}
   One may alternatively use the
Dynkin-Helffer-Sj\"ostrand formula as in the Appendix, and both definitions
are known to coincide on $L^2(\Omega)$. However, the
Dynkin-Helffer-Sj\"ostrand formula seems to be restricted to defining
$m(-\Delta_D)$ for functions $m$ which exhibit slightly more decay
than required in the next theorem, at least if one proceeds as exposed
in the Appendix.
\end{remark}
\begin{theorem}[\cite{Alex1}]
\label{mikhlin}
  Let $m\in C^{N}(\R^+)$, $N\in\N$ and $N\geq n/2+1$, such that
  \begin{equation}
    \label{eq:mm}
    \sup_{\xi,k\leq N}|\xi \partial_\xi^k m(\xi)|<+\infty.
  \end{equation}
Then the operator defined by \eqref{eq:sM} extends to a continuous
operator on $L^p(\Omega)$, and sends $L^1(\Omega)$ to weak
$L^1(\Omega)$.
\end{theorem}
In order to use the argument of \cite{Alex1}, we need the
Gaussian upper bound on the heat kernel, which is provided in our case
by \cite{DavHeat}. Once we have Theorem \ref{mikhlin}, all we need to
do to prove Theorem \ref{theoremtriebel} is to follow Stein's
classical proof from \cite{Stein}\footnote{we thank Hart Smith for bringing
this to our attention}, and we
recall it briefly for the convenience of the reader. Let us introduce
the Rademacher functions, which are defined as follows:
\begin{itemize}
\item the function $r_0(t)$ is defined by $r_0(t)=1$ on $[0,1/2]$ and
  $r_0(t)=-1$ on $(1/2,1)$, and then extended to $\R$ by periodicity;
\item for $m\in \N\setminus\{0\}$, $r_m(t)=r_0(2^m t)$.
\end{itemize}
Their importance is outlined by the following inequalities (see the
Appendix in \cite{Stein}),
\begin{equation}
  \label{eq:kin}
  c_p \| \sum_m a_m r_m(t) \|_{L^p_t} \leq (\sum_m |a_m|^2)^\frac 1 2
\leq C_p \| \sum_m a_m r_m(t) \|_{L^p_t}.
\end{equation}
Now, define
$$
m^\pm(t,\xi)=\sum_{j=0}^{+\infty} r_j(t)\Psi_{\pm j}(\xi),
$$
where $\Psi_j$ was defined in the introduction. A straightforward
computation proves that the bound \eqref{eq:mm} holds for
$m^\pm(t,\xi)$. Therefore, 
$$
\| m^\pm(t,-\Delta_D)f\|_{L^p(\Omega)} \lesssim \|f\|_{L^p(\Omega)};
$$
integrating in time over $[0,1]$, exchanging space and time norms, and using \eqref{eq:kin}, 
$$
\|m^\pm(t,-\Delta_D) f\|_{L^p(\Omega) L^2(0,1)}\approx
\|\big (\sum_{j=0}^{\pm\infty} |\Psi(-2^{-2j}\Delta_D) f|^2 \big
)^\frac 1 2\|_{L^p(\Omega)} \lesssim \|f\|_{L^p(\Omega)}.
$$
This proves one side of the equivalence in \eqref{ib}: the other side
follows from duality, once we see the above estimate as an estimate from
$L^p(\Omega)$ to $L^p(\Omega;l^2)$, which maps $f$ to
$(\Psi(-2^{-2j}\Delta_D) f)_{j\in \Z}$.
\section{Heat flow estimates}
In order to prove Proposition \ref{propgrad} we need the following
lemma.
\begin{lemma}\label{lemest1}
For all $1\leq p\leq +\infty$, we have
\begin{equation}
  \label{eq:limzero}
\|S(t)f\|_{L^{p}(\Omega)}\rightarrow_{t\rightarrow\infty} 0,
\end{equation}
\begin{equation}
  \label{eq:bornesup}
 \sup_{t\geq 0}\|S(t)f\|_{L^{p}(\Omega)}\lesssim\|f\|_{L^{p}(\Omega)}\,.
\end{equation}
Moreover,
\begin{equation}
  \label{eq:bornemax}
\|\sup_{t\geq 0}|S(t)f|
\|_{L^{p}(\Omega)}\lesssim\|f\|_{L^{p}(\Omega)}\lesssim \|f\|_{L^{p}(\Omega)}\,.
\end{equation}
\end{lemma}
\emph{Proof:} The estimate \eqref{eq:bornesup} clearly follows from
\eqref{eq:bornemax}, which in turn is a direct consequence of the
Gaussian nature of the Dirichlet heat kernel, see \cite{DavHeat}. The
same Gaussian estimate implies \eqref{eq:limzero}. However we do not
need such a strong fact to prove  \eqref{eq:bornesup}, which will
follow from the next computation as well (see
\eqref{eq:base}) when $1<p<+\infty$. Estimate \eqref{eq:limzero} can also be obtained
through elementary arguments. We defer such a proof to the end of the section.

\subsection{Proof of Theorem \ref{propgrad}}

If $p=2$ the proof is nothing more than the energy inequality,
combined with \eqref{eq:limzero}. In fact, for $p=2$, we have equality
in \eqref{eq:quadheat} with $C_2=2$. We now take $p=2m$ where $m\geq
2$. Multiplying equation \eqref{heateq} by  $\bar{u}|u|^{p-1}$ and
taking the integral over $\Omega$ and $[0,T]$, $T>0$ yields, taking
advantage of the Dirichlet boundary condition,
\begin{multline}
  \label{eq:base}
  \frac{1}{p}\int_{0}^{T}\partial_{t}\|u\|^{p}_{L^{p}(\Omega)}dt+\int_{0}^{T}\int_{\Omega}|\nabla
u|^{2}|u|^{p-2}dxdt+\\{}+\frac{(p-2)} 2 \int_{0}^{T}\int_{\Omega}(\nabla
(|u|^2))^{2}|u|^{p-4}dxdt=0,
\end{multline}
from which we can estimate either $\|u\|^{p}_{L^{p}(\Omega)}(T)\leq
\|f\|^{p}_{L^{p}(\Omega)}$ (which is \eqref{eq:bornesup}) or
\[
\|f\|^{p}_{L^{p}(\Omega)}\leq\|u\|^{p}_{L^{p}(\Omega)}(T)+p(p-1)\int_{0}^{T}\int_{\Omega}|\nabla u|^{2}|u|^{p-2}dxdt.
\]
Letting $T$ go to infinity and using \eqref{eq:limzero} from Lemma \ref{lemest1} and H\"older inequality we find
\[
\|f\|^{p}_{L^{p}(\Omega)}\leq
p(p-1)\Big(\int_{\Omega}\Big(\int_{0}^{\infty} |\nabla
u|^{2}\,dt\Big)^{\frac p 2} \,dx\Big)^{\frac 2
  p}\Big(\int_{\Omega}(\sup_t |u|^{p-2})^{\frac p{p-2}}dx\Big)^{\frac {p-2}p}\,.
\]
The proof follows using again Lemma \ref{lemest1}, as
\[
\|f\|^{p}_{L^{p}(\Omega)}\leq C_{p}\|\Big(\int_{0}^{\infty}|\nabla
u|^{2}dt\Big)^\frac 1 2\|_{L^{p}(\Omega)}\left(\|\sup_{t\geq 0}|u|\|_{L^{p}(\Omega)}\right)^{p-2}.
\] 
Note that we may prove the weaker part, \eqref{eq:quadbheat}, without
assuming the maximal in time bound, by reversing the order of
integration in our argument. This would keep the argument for heat
square functions essentially
self-contained, without any need for Gaussian bounds on the heat kernel.
\begin{remark}
  We do not claim novelty here: our argument follows closely (a dual version
  of) the proof of a classical square function bound for the Poisson kernel in
  the whole space, see \cite{Stein}.
\end{remark}
We have proved one side of the equivalence in \eqref{eq:quadheat}
involving the $Q_t$ square function, in the range $2\leq p<+\infty$; we now prove the other side, by
duality. Let $\phi \in L^q(\Omega)$, with $1/q=1-2/p$, and consider
$$
I=\int_\Omega \big( \int_0^{+\infty}|\nabla u|^2\,dt\big) \phi(x)\,dx.
$$
Without any loss of generality, we may assume $\phi\geq 0$. On the
other hand, let $v=|\nabla u|^2$, then
$$
\partial_t v-\Delta v=-2|\nabla^2 u|^2,
$$
and one checks easily that $\partial_n v=0$ on $\partial\Omega$. Let
$S_n(t)$ be the solution to the heat equation on $\Omega$ with Neumann
boundary condition, by comparing $v$ and $S_n(t/2) v(t/2)$ (formally,
take the difference, multiply by the positive part and integrate by parts) we have
$$
0\leq v\leq S_n(t/2)v(t/2)=S_n(t/2)|\nabla u(t/2)|^2,
$$
and therefore
\begin{equation}
  \label{eq:dualheat}
  I\leq 2 \int_\Omega \int_0^{+\infty}|\nabla u|^2
\,S_n(t) \phi\,dxdt.
\end{equation}
Now, we also have
$$
\partial_t u^2-\Delta_D u^2=-2|\nabla u|^2,
$$
and therefore
$$
I\leq  -\int_\Omega \int_0^{+\infty} (\partial_t -\Delta)( u^2)
S_n(t) \phi\,dxdt.
$$
From $(\partial_t -\Delta)( u^2 S_n(t) \phi)=-2 \nabla(u^2)\cdot
\nabla S_n(t) \phi$, we get
$$
I\leq \int_\Omega 4 \sup_t |u| \big(\int_0^{+\infty} |Q_t u|^2 \,
\frac{dt}{t}\big)^\frac  1 2 \big(\int_0^{+\infty} |\nabla S_n(t) \phi|^2 \,
{dt}\big)^\frac  1 2.
$$
The bound we already proved with $Q_t$ can easily be reproduced with
$S(t)$ replaced by $S_n(t)$, and therefore, provided $q\leq 2$, we may
use the dual bound on the square function of $\phi\in L^q(\Omega)$ and
conclude by H\"older, using \eqref{eq:bornemax} on the first factor. The condition on $q$ translates into $p\geq 4$,
and the remaining $2<p<4$ are handled by interpolation.
\begin{remark}
  Actually, we may directly bound $S_n(t)$ by a Gaussian in
  \eqref{eq:dualheat}, extend $\phi$ by $0$ outside $\Omega$, and use
  the heat square function bounds in $\R^n$. This provides a direct
  argument, irrespective of the value of $p$.
\end{remark}
It remains to prove the equivalence between the $Q_t$ square function
and the $\mathrm{\bf Q}_t$ square function. For this, we repeat the
duality argument but we replace $|\nabla u|^2$ by $t |\partial_t
u|^2$. Notice that $\partial_t u$ is also a solution to the heat
equation with Dirichlet boundary condition, and if $w=|\partial_t u|^2$,
$$
(\partial_t-\Delta) w=-2|\nabla \partial_t u|^2.
$$
Therefore, comparing $w$ and $S(t/2)w(t/2)$,
$$
0\leq |\partial_t u|^2 \leq S(t/2)|\partial_t u(t/2)|^2,
$$
and
$$
J=\int_\Omega \int_0^{+\infty}  |\partial_t u|^2
 t \phi\,dxdt \leq 2  \int_\Omega \int_0^{+\infty}  |\partial_t u|^2
 t S(t) \phi\,dxdt.
$$
Now,
\begin{align*}
  J & \leq 2\int_\Omega \int_0^{+\infty} t \partial_t u\Delta u
  S(t) \phi\,dxdt\\
 & \leq -\int_\Omega \int_0^{+\infty} t \partial_t |\nabla u|^2
  S(t) \phi\,dxdt-2 \int_\Omega \int_0^{+\infty} t \partial_t u \nabla u
  \nabla S(t) \phi\,dxdt\\
 & \leq \int_\Omega \int_0^{+\infty} |\nabla u|^2
  (1+t\partial_t) S(t) \phi\,dxdt+2 \int_\Omega \int_0^{+\infty}
  |\mathrm{\bf Q}_t u Q_t u
  Q_t \phi|\,\frac{dt} t \, dx
\end{align*}
from which we can easily conclude by H\"older (using Lemma
\ref{lemest} to bound $t\partial_t S(t)\phi$). Duality takes care of
the reverse bound, and this concludes the proof of Theorem
\ref{propgrad}, except for the equivalence between the $Q_t$ and
$\mathrm{\bf Q}_t$ Besov norms in \eqref{eq:quadbheat}; we defer this
to the end of  the next subsection.

Notice that, at this point, we proved Theorem \ref{theorembesov}, but
with the $\Psi$ operator replaced by the gradient heat kernel and the
discrete parameter $2^{-2j}$ by the continuous parameter $t$. The rest
of this section is devoted to proving the equivalence between the Besov
norms which are defined by the heat kernel or the spectral localization.
\begin{lemma}\label{lemequiv}
Let $1\leq p\leq +\infty$. We have the following equivalence between dyadic and continuous
versions of the Besov norm: 
\[
\frac 3 4 \sum_{k\in\Z}\|Q_{2^{-2k}}f\|^{2}_{L^{p}(\Omega)}\leq \int_{0}^{\infty}\|Q_tf\|^{2}_{L^{p}(\Omega)}\frac{dt}{t}\leq 3 \sum_{k\in\Z}\|Q_{2^{-2k}}f\|^{2}_{L^{p}(\Omega)}.
\]
\end{lemma}
This follows at once from factoring the semi-group: for
$2^{-2j}\leq t\leq 2^{-2(j-1)}$, write $S(t)=S(t-2^{-2j})S(2^{-2j})$
and use \eqref{eq:bornesup}. We now turn to the direct proof of Theorem
\ref{theorembesov} from the heat flow version. Let $\Psi\in
C^{\infty}_{0}(\mathbb{R}^{*})$ satisfying \eqref{unu} and denote
$\Delta_{j}f\eqdef \Psi(2^{-2j}\Delta_{D})f$, where
$\Psi(2^{-2j}\Delta_{D})f$ is given by the
Dynkin-Helffer-Sj\"{o}strand formula (see the Appendix, \eqref{integralpsi}). 
From Proposition \ref{propgrad} and Lemma \ref{lemequiv} we have
\begin{equation}\label{q}
\|f\|_{L^{p}(\Omega)}\leq 3 C_{p}\Big(\sum_{k\in\Z}\|Q_{2^{-2k}}f\|^{2}_{L^{p}(\Omega)}\Big)^{1/2}
\end{equation}
and we will show that \eqref{q} implies \eqref{ib}: it suffices to
prove the following almost orthogonality property between localization
operators $\Delta_j$ and $Q_{2^{-2k}}$:
\begin{equation}\label{estkj}
\forall k,j\in\Z,\quad\quad\|Q_{2^{-2k}}\Delta_{j}f\|_{L^{p}(\Omega)}\lesssim 2^{-|j-k|}\|\Delta_{j}f\|_{L^{p}(\Omega)}.
\end{equation}
Then, from $(2^{-|j-k|})_{k}\in l^{1}$ and
$(\|\Delta_{j}f\|_{L^{p}(\Omega)})_{j}\in l^{2}$ we estimate
\begin{equation}\label{estqd}
\sum_{k\in\Z}\|Q_{2^{-2k}}f\|^{2}_{L^{p}(\Omega)}=\sum_{k\in\Z}\|\sum_{j\in\Z}Q_{2^{-2k}}\Delta_{j}f\|^{2}_{L^{p}(\Omega)}
\end{equation}
as an $l^{1}*l^{2}$ convolution and conclude using Lemma \ref{lemest}.
It remains to show \eqref{estkj}:
\begin{itemize}
\item for $k<j$ we write
  \begin{multline*}
    Q_{2^{-2k}}\Delta_{j}f=2^{3/2}2^{-2(j-k)}\Big(2^{-(2k+1)/2}\nabla
S(2^{-(2k+1)})\Big)\\\Big(2^{-(2k+1)}\Delta_{D}S(2^{-(2k+1)})\Big)\breve{\Psi}(-2^{-2j}\Delta_{D}){\Psi}(-2^{-2j}\Delta_{D})f,
  \end{multline*}
where we set $\breve{\Psi}(\lambda)\eqdef \frac{1}{\lambda}\tilde
\Psi(\lambda)$, and $\tilde{\Psi}\in C_0^\infty$, $\tilde{\Psi}=1$ on
$\mathrm{supp} \Psi$. By Lemma \ref{lemest}, the operators
$Q_{2^{-(2k+1)}}=2^{-(2k+1)/2}\nabla S(2^{-(2k+1)})$ and
$\mathrm{\bf Q}_{2^{-(2k+1)}}=2^{-(2k+1)}\Delta_{D}S(2^{-(2k+1)})$ are
bounded on $L^{p}(\Omega)$ and we obtain \eqref{estkj} using Corollary
\ref{cor} for $\breve{\Psi}$.
\item for $k\geq j$ we set $\Psi_{1}(\xi)=\tilde{\Psi}(\xi)\exp(\xi)$,
  $\Psi_{2}(\xi)=\Psi(\xi)$, and we use again Lemma \ref{lemsupp} to write (slightly abusing the
notation as $2^{-2k}-2^{-2j}<0$)
\begin{equation}
  \label{eq:backheat}
  S(2^{-2k}-2^{-2j})\Delta_{j}f=S(2^{-2k}){\Psi_1}(-2^{-2j}\Delta_{D}) \Psi_2(-2^{-2j}\Delta_{D})f.
\end{equation}
Then
\[
Q_{2^{-2k}}\Delta_{j}f=2^{-(k-j)}\Big(2^{-j}\nabla S(2^{-2j})\Big)\Big(S(2^{-2k}-2^{-2j})\Delta_{j}f\Big),
\]
 and using again Lemma \ref{lemest} we see that the operator $2^{-j}\nabla
S(2^{-2j})$ is bounded while the remaining operator \eqref{eq:backheat}
is bounded by Corollary \ref{cor}. This ends the proof.
\end{itemize}
\begin{remark}\label{eqBesov}
  One may prove a similar bound with $Q_{2^{-2k}}$ and $\Delta_j$
  reversed, either directly or by duality. Hence Besov norms
  based on $\Delta_j$ or $Q_{2^{-2k}}$ are equivalent.
\end{remark}

\subsection{Proof of Proposition \ref{lemest}}
For $\mathrm{\bf
  Q}_t$, boundedness  on all $L^p$ spaces, including $p=1,+\infty$, follows once again from a
Gaussian upper bound on $\partial_t S(t)$ (see \cite{DaviesJOT} or
\cite{DaviesNGA}). However the subsequent Gaussian bound on the
gradient $\nabla_x S(t)$ in \cite{DaviesJOT} is a direct consequence
of the Li-Yau inequality, which holds only inside convex domains. We
were unable to find a reference which would provide the desired bound
for $Q_t$ in the context of the exterior domain. Therefore we provide an elementary detailed
proof for $Q_{t}$. Furthermore, we only deal with $1<p<2$ or powers of two, $p=2^m$,
$m\in\N^{*}$: complex interpolation takes care of remaining values of
$p$, though one could adapt the following argument to generic values $p>2$,
at the expense of lengthier computations.

Set $v(x,t)=(v_{1},..,v_{n})(x,t):=Q_{t}f=t^{1/2}\nabla u(x,t)$ and assume without loss of
generality that $v_{j}$ are real: we multiply the equation satisfied by $v$
by $v|v|^{p-2}$, where $|v|^{2}=\sum_{j=1}^{n}v_{j}^{2}$, and integrate over $\Omega$,
 \begin{multline}\label{eqv}
\partial_t \left( \frac{1}{p}\|v\|^{p}_{L^{p}(\Omega)}\right)-\sum_{j=1}^{n}\int_{\partial\Omega}((\overrightarrow{\nu}\cdot
\nabla) v_{j})\cdot v_{j} |v|^{p-2}d\sigma +  \\ + \int_{\Omega}|\nabla v|^{2}|v|^{p-2}\,dx+\frac{(p-2)}{2}
\int_{\Omega}\nabla (|v|^{2})|v|^{p-4}\,dx= \frac{1}{2t}\|v\|^{p}_{L^{p}(\Omega)},
\end{multline}
where $\overrightarrow{\nu}$ is the outgoing unit normal vector to $\partial\Omega$ and $d\sigma$ is the surface measure on $\partial\Omega$. We claim that the second term in the left hand side vanishes: in fact we write
 \begin{multline}
\sum_{j=1}^{n}\int_{\partial\Omega}(\overrightarrow{\nu}\cdot \nabla v_{j}) \cdot v_{j} |v|^{p-2}d\sigma=\\=\frac{t^{p/2}}{2}\int_{\partial\Omega}\partial_{\nu}(|\partial_{\nu}u|^{2}+|\nabla_{tang}u|^{2})(|\partial_{\nu}u|^{2}+|\nabla_{tang}u|^{2})^{(p-2)/2}d\sigma,
 \end{multline}
and from $u|_{\partial\Omega}=0$ the time and tangential
derivative $(\partial_t,\nabla_{\text{tang}})u|_{\partial\Omega}$
vanishes; furthermore, using the equation, $\partial^2_{\nu}u=0$ on
$\partial\Omega$. 
\begin{remark}
Notice that while this term does not vanish with Neumann boundary
conditions, it will be a lower order term (like $|\nabla u|^2$ on
$\partial\Omega$) which can be controled by the trace theorem.
\end{remark}

Now, if $1<p<2$, multiply by $\|v\|^{2-p}_{L^p(\Omega)}$ and integrate
over $[0,T]$,
\begin{equation*}
  \|v\|^2_{L^p(\Omega)}(T) \lesssim \int_0^T \| Q_t
  f\|^2_{L^p(\Omega)}\,\frac{dt}{t}\lesssim \|f\|^2_p,
\end{equation*}
where the last inequality is the dual of \eqref{eq:quadheat}. Hence we
are done with $1<p<2$.
\begin{remark}
  We ignored the issue of $v$ vanishing in the third term in
  \eqref{eqv}. This is easily fixed by replacing $|v|^{p-2}$ by
  $(\sqrt{\varepsilon+|v|^2})^{p-2}$ and proceeding with the exact same
  computation. Then let $\varepsilon$ go to $0$ after dropping the
  positive term on the left handside of \eqref{eqv}.
\end{remark}
Now let $p=2^m$ with $m\geq 1$: we proceed directly by integrating
\eqref{eqv} over $[0,T]$, to get
 \begin{multline}\label{eqvbis}
 \frac{1}{p}\|v\|^{p}_{L^{p}(\Omega)}(T)+\int_0^T \int_{\Omega}|\nabla v|^{2}|v|^{p-2}\,dxdt +\\ + \frac{(p-2)}{2}
\int_0^T \int_{\Omega}|\nabla (|v|^{2})|^{2}|v|^{p-4}\,dxdt =\int_0^T \frac{1}{2t}\|v\|^{p}_{L^{p}(\Omega)}\,dt.
\end{multline}
On the other hand (recall \eqref{eq:base}),
\begin{equation}\label{equ}
\frac{1}{p}\|u\|^{p}_{L^{p}(\Omega)}(T)+(p-1)\int_0^T
\int_{\Omega}|\nabla u|^{2}|u|^{p-2}\,dxdt=\frac 1 p \|f\|^p_{L^p(\Omega)}.
\end{equation}
If $p=2$ the estimates are trivial since from \eqref{eqvbis}, \eqref{equ},
\[
\frac 1 2 \|v\|^{2}_{L^{2}(\Omega)}(T)\leq
\int_{0}^{T}\frac{1}{2t}\|v\|^{2}_{L^{2}(\Omega)}dt =\frac 1 2
\int_0^T \|\nabla u\|^2_{L^2(\Omega)}\,dt\leq \frac 1 4 
\|f\|^{2}_{L^{2}(\Omega)}.
\]
Now, let $p\geq 4$; for convenience, denote by $J$ the second integral in the left hand-side
of \eqref{eqvbis} (notice that the third integral is bounded from above by $J$), hence
$$
J=
\int_0^{T}\int_{\Omega} |\nabla^2 u|^2 |\nabla u|^{p-2}
t^{\frac p 2}\,dxdt=\int_{0}^{T}\int_{\Omega}(\sum_{i,j}|\partial^{2}_{i,j}u|^{2})(\sum_{j}|\partial_{j}u|^{2})^{\frac{(p-2)}{2}}t^{\frac{p}{2}}\,dx dt,
$$
and set
\begin{equation}\label{defIk}
I_k=\int_{0}^{T}\int_{\Omega}|\nabla u|^{2k}|u|^{p-2k}
t^{k -1}\,dxdt \,\,\text{where}\,\, 2\leq 2k\leq p.
\end{equation}
For our purposes, it suffices to estimate the 
right hand-side of \eqref{eqvbis}, which rewrites
\begin{equation}
  \label{eq:I}
\frac 1 2  \int_0^{T} t^{\frac p 2 -1} \|\nabla
u\|^p_{L^p(\Omega)}\,dt=\frac 1 2 I_{\frac p 2}.
\end{equation}
Integrate by parts the inner (space) integral in $I_k$, the boundary
term vanishes and collecting terms,
\begin{multline}\label{pegal4}
\int_\Omega \nabla u \nabla u |\nabla
u|^{2(k-1)}|u|^{p-2k}\,dx\leq \frac{(2k-1)}{(p-2k+1)}\int_{\Omega}|\nabla^{2}u||\nabla
u|^{2k-2}|u|^{p-2k+1}dx.
\end{multline}
By Cauchy-Schwarz the integral in the right hand side of \eqref{pegal4} is bounded by
\[
\Big(\int_{\Omega}|\nabla^{2}u|^{2}|\nabla
u|^{p-2}dx\Big)^{1/2}\Big(\int_{\Omega}|\nabla
u|^{4k-4-(p-2)}|u|^{2p+2-4k}dx\Big)^{1/2},
\]
therefore for $k\geq \frac{p}{4}+1$ we have
\[
I_k\lesssim \frac{(2k-1)}{(p-2k+1)} J^{\frac 1 2} I_{2k-\frac p 2-1}^\frac 1 2.
\]
We aim at controlling $I_{m}$ by $J^{1-\eta}I_1^\eta$, for some
$\eta>0$ which depends on $m$ (notice that when $p=4$, which is $m=2$, we are already
done, using $k=2$ !). Set $k=\frac p 2-(2^j-1)$ with $j\leq m-2$,
$$
I_{2^{m-1}-(2^j-1)} \leq\frac{(2^{m}-(2^{j+1}-1))}{(2^{j+1}-1)} J^\frac 1 2 I_{2^{m-1}-(2^{j+1}-1)}^\frac
1 2,
$$
and iterating $m-2$ times, we finally control $I_\frac p 2$ by
$J^{1-\eta}I_1^\eta$, which proves that $Q_t$ is bounded on
$L^p(\Omega)$.

We now proceed to obtain boundedness of $\mathrm{\bf Q}_{t}$ on
$L^{p}(\Omega)$ from the $Q_t$ bound; this is worse than
using the Gaussian properties of its kernel, as the constants blow up
when $p\rightarrow 1,+\infty$. It is, however, quite simple. By
duality $Q^\star_t$ is bounded on $L^p(\Omega)$, and
$$
\mathrm{\bf Q}_t=t\partial_t S(t)=tS(\frac t 2)\Delta S(\frac t
2)=2 \sqrt \frac t 2 S(\frac t 2)\nabla\cdot \sqrt \frac t 2 \nabla
S(\frac t 2)=2
Q^\star_{\frac t 2} Q_{\frac t 2},
$$
and we are done with Lemma \ref{lemest}.

From the previous decomposition, we also obtain 
$$
\|\mathrm{\bf Q}_t f\|_{L^p(\Omega)}\lesssim \|Q_t
f\|_{L^p(\Omega)},
$$
which implies that any Besov norm defined with $\mathrm{\bf Q}_t$ is
bounded by the corresponding norm for $Q_t$. The reverse bound is true
as well, though slightly more involved. We provide the proof for
completeness. Consider $f,h\in C^\infty_0(\Omega)$ and $\langle
f,g\rangle=\int_\Omega fg$. Then
\begin{align*}
  \langle
f,g\rangle & = - \int_0^{+\infty} \langle \partial_t S(t) f,h\rangle
\,dt = -2 \int_0^{+\infty} \langle \partial_t S(t) f,S(t) h\rangle
\,dt\\
 & = 2 \int_{t<s} \langle \partial_t S(t) f,\partial_s S(s) h\rangle
\,dtds = 4 \int_0^{+\infty} \langle \nabla S(s)\partial_t S(t) f,\nabla
S(s) h\rangle\,dtds\\
 &\lesssim \int_s \left\|\int_0^s \nabla S(t) \partial_s S(s) f
 \,dt\right \|_p \|\nabla S(s) h\|_{p'}\,ds\lesssim \int_s \sqrt s\| \partial_s S(s) f
  \|_p \|\nabla S(s) h\|_{p'}\,ds
\end{align*}
where we used our bound on $\sqrt t \nabla S(t)$ at fixed $t$. Then
$$
\langle f,h\rangle  \lesssim \int_s \| \mathrm{\bf Q}_s f
\|_p \|Q_s h\|_{p'}\,\frac{ds}s
$$
from which we are done by H\"older.

\subsection{Proof of Proposition \ref{lemestvect}}
Let us consider now the vector valued case $u=(u_{l})_{l\in\{1,..,N\}}$ for $N\geq
2$, where each $u_{l}$ solves \eqref{heateq} with Dirichlet condition
and initial data $f_{l}$. For the sake of simplicity we consider only
real valued $u_{l}$, and write 
$$
|u|^{2}=\sum_{l=1}^{N}u_{l}^{2},\,\,\, |\nabla
u_{l}|^{2}=\sum_{j=1}^{n}(\partial_{j}u_{l})^{2},\,\,\,|\nabla
u|^{2}=\sum_{j=1}^{n}\sum_{l=1}^{N}(\partial_{j}u_{l})^{2}
$$
Notice that $n$ is the spatial dimension and is fixed through the
argument: hence all constants may depend implicitely on $n$, while $N$
is the dimension of $H$. For $p=1,+\infty$, the boundedness of
$\mathrm{\bf Q}_t$ follows
from the Gaussian character of the time derivative heat kernel, which
is diagonal on $H$.

We proceed with $Q_t$. Multiplying the
equation satisfied by $u_{l}$ by $u_{l}|u|^{p-2}$, integrating over
$\Omega$ and  summing up we immediately get \eqref{eq:base}. We now
proceed to obtain bounds for $v(x,t)=(v_{l}(x,t))_{l}$,
where $v_{l}(x,t)=t^{1/2}\nabla u_{l}(x,t)$. 
Multiplying the equation satisfied by $v_{l}$ by $v_{l}|v|^{p-2}$
where $|v|^{2}=t|\nabla u|^{2}$, summing up over $l$ and taking the
integral over $\Omega$ yields
\begin{multline}
\frac{1}{p}\|v\|^{p}_{L^{p}(\Omega)}(T)+\sum_{l=1}^{N}\sum_{j=1}^{n}\int_{0}^{T}\int_{\Omega}|\nabla(\partial_{j}u_{l})|^{2}|\nabla
u|^{p-2}\,dxdt+\\+\frac{(p-2)}{4}\int_{0}^{T}\int_{\Omega}|\nabla|\nabla
u|^{2}|^{2} |\nabla u|^{p-4}t^{p/2}\,dxdt=\int_{0}^{T}\|\nabla
u\|^{p}_{L^{p}(\Omega)}t^{p/2-1}\, dt=\frac{1}{2}I_{\frac{p}{2}},
\end{multline}
where
$|\nabla(\partial_{j}u_{l})|^{2}=\sum_{i=1}^{n}(\partial^{2}_{i,j}u_{l})^{2}$,
$|\nabla
u|^{2}=\sum_{l=1}^{N}\sum_{j=1}^{n}(\partial_{j}u_{l})^{2}$. Notice
again that the boundary term vanishes. Denote the last two integrals
in the left hand side
by $J_{1}$, $J_{2}$.
Like before, we perform integrations by parts in $I_{k}$ defined in \eqref{defIk} to obtain
\begin{multline}\label{term4line2}
\int_{\Omega}|\nabla u|^{2k}|u|^{p-2k}\,dx=-\sum_{l=1}^{N}\int_{\Omega}u_{l}\Delta u_{l}|\nabla u|^{2(k-1)}|u|^{p-2k}\,dx-\\-(k-1)\sum_{l=1}^{N}\int_{\Omega}u_{l}\nabla u_{l}\nabla(|\nabla u|^{2})|\nabla u|^{2(k-2)}|u|^{p-2k}\,dx-\\-(p-2k)\sum_{i=1}^{n}\int_{\Omega}(\sum_{l=1}^{N}\partial_{i}u_{l}u_{l})^{2}|u|^{p-2k-2}\,dx.
\end{multline}
For $k\geq \frac{p}{4}+1$ we estimate the first term in the right hand side of \eqref{term4line2} by
\begin{multline*}
  \int_{\Omega} (\sum_{l=1}^{N}u_{l}^{2})^{1/2}(\sum_{l=1}^{N}(\Delta u_{l})^{2})^{1/2}|\nabla u|^{2(k-1)}|u|^{p-2k}\,dx\leq\\
(\sum_{l=1}^{N}\sum_{j=1}^{n}\int_{\Omega}\int_{\Omega}|\nabla(\partial_{j}u_{l})|^{2}|\nabla u|^{p-2}\,dx)^{1/2}(\int_{\Omega}|\nabla u|^{4k-p-2}|u|^{2p-4k+2}\,dx)^{1/2},
\end{multline*}
and the second term in the right hand side of \eqref{term4line2} by
\[
(k-1)(\int_{\Omega}\sum_{i=1}^{n}(\partial_{i}(|\nabla u|^{2}))^{2}|\nabla u|^{p-4}\,dx)^{1/2}(\int_{\Omega}|\nabla u|^{4k-p-2}|u|^{2p-4k+2}\,dx)^{1/2},
\]
where we used that 
\[
\sum_{l=1}^{N}u_{l}\nabla u_{l}\nabla(|\nabla u|^{2})\leq \sum_{i=1}^{n}(\sum_{l=1}^{N}u_{l}^{2})^{1/2}(\sum_{l=1}^{N}(\partial_{i}u_{l})^{2})^{1/2}|\partial_{i}(|\nabla u|^{2})|\lesssim |u||\nabla u||\nabla(|\nabla u|^{2})|.
\]
Since the last term in \eqref{term4line2} is negative, while the
quantity we want to estimate is positive we obtain from the last
inequalities
\begin{equation}\label{inegIk}
\int_{0}^{T}\int_{\Omega}|\nabla u|^{2k}|u|^{p-2k}t^{k-1}\,dxdt\lesssim (J_{1}^{1/2}+J_{2}^{1/2})I_{2k-\frac{p}{2}-1}\lesssim (J_{1}+J_{2})^{1/2}I_{2k-\frac{p}{2}-1}.
\end{equation}
From now on we proceed exactly like in the scalar case iterating
sufficiently many times to obtain the desired result, since we control
$I_{p/2}$ which is the RHS term of \eqref{term4line2} using
\eqref{inegIk}.

\subsection{A simple argument for \eqref{eq:limzero}}
We now return to the first estimate in Lemma \ref{lemest1}: while we
only deal with $p=2$, there is nothing
specific to the $L^2$ case in what follows. Let $\chi$ be a smooth
cut-off near the boundary $\partial \Omega$. Then $v=(1-\chi)u$ solves
the heat equation in the whole space, with source term $[\chi,\Delta]
u$:
$$
(1-\chi)u =S_0(t) (1-\chi)u_0+\int_0^t S_0(t-s) [\chi,\Delta] u(s)
\,ds,
$$
where $S_0$ is the free heat semi-group. We have, taking advantage of
the localization near the boundary,
$$
\|[\chi,\Delta] u\|_{L^2_t (L^{\frac {2n}{n+2}})}\lesssim C(\chi,\chi')
  \|\nabla u\|_{L^2_t (L^{2})}<+\infty,
$$
by the energy inequality \eqref{eq:base}. The integral equation
on $(1-\chi)u$ features $S_0$ for which we have trivial Gaussian
estimates, and both the homogeneous and inhomogeneous terms are
$C_t(L^2)$ and go to zero as time goes to $+\infty$. On the other
hand, by Poincar\'e inequality (or Sobolev),
$$
\int_0^t \|\chi u\|^2_2 \,ds\lesssim \int_0^t \|\nabla u\|^2_2 \,ds,
$$
which ensures that $\|\chi u\|_2$ goes to zero as well at $t=+\infty$.

\section*{Acknowledgments}
The authors would like to thank Hart Smith for pointing out the
relevance of \cite{Alex1,Alex2,DOS} in this context, Pascal Auscher and Francis Nier for
entertaining discussions, not to mention providing material which
greatly improved content, and Nikolay Tzvetkov for helful remarks on an early draft. Part of this work was conducted while the second author was visiting
the Mittag-Leffler institute, which he is grateful to for its
hospitality. Both authors were partially supported by the A.N.R. grant ``Equa-disp''.

\section*{Appendix: functional calculus}
We start by recalling the Dynkin-Helffer-Sj\"{o}strand formula
(\cite{Dynkin,HelSj}) and refer to the appendix of \cite{francis} for
a nice presentation of the use of almost-analytic extensions in the
context of functional calculus. In what follows we will also rely on Davies'presentation (\cite{Dav89}) from which we will use a couple of
useful lemma.
\begin{definition}(see \cite[Lemma A.1]{francis})
Let $\Psi\in C^{\infty}_{0}(\mathbb{R})$, possibly complex valued. We
assume that there exists $\tilde{\Psi}\in C^{\infty}_{0}(\mathbb{C})$
such that $|\bar{\partial}\tilde{\Psi}(z)|\leq C|\Im z|$ and
$\tilde{\Psi}|_{\mathbb{R}}=\Psi$. Then we have (as a bounded operator
in $L^2(\Omega)$)
\begin{equation}\label{integralpsi}
\Psi(-h^{2}\Delta_{D})=\frac{i}{2\pi}\int_{\mathbb{C}}\bar{\partial}\tilde{\Psi}(z)(z+h^{2}\Delta_{D})^{-1}d\bar{z} \wedge dz.
\end{equation}
\end{definition}
The next result ensures the existence of $\tilde{\Psi}$ in the
previous definition ( see \cite[Lemma A.2]{francis} and \cite{trev},
where it is linked with Hadamard's problem of finding a smooth function
with prescribed derivatives at a given point):
\begin{lemma}
If $\Psi$ belongs to $C^{\infty}_{0}(\mathbb{R})$ there exists $\tilde{\Psi}\in C^{\infty}_{0}(\mathbb{C})$ such that $\tilde{\Psi}|_{\mathbb{R}}=\Psi$ and 
\begin{equation}\label{estimaextpsi}
|\bar{\partial}\tilde{\Psi}(z)|\leq C_{N,\Psi}|\Im z|^{N},\quad \forall z\in \mathbb{C}, \quad \forall N\in\mathbb{N}.
\end{equation}
Moreover, if $\Psi$ belongs to a bounded subset of
$C^{\infty}_{0}(\mathbb{R})$ (elements of $\mathcal{B}$ are
supported in a given compact subset of $\mathbb{R}$ with uniform bounds), then the mapping $\mathcal{B}\ni \Psi\rightarrow \tilde{\Psi}\in
C^{\infty}_{0}(\mathbb{C})$ is continuous and $C_{N,\Psi}$ can be
chosen uniformly w.r.t $\Psi\in \mathcal{B}$. 
\end{lemma}
\begin{remark}
Estimate \eqref{estimaextpsi} simply means that $\bar{\partial}\tilde{\Psi}(z)$ vanishes at any order on the real axis. Precisely, if $z=x+iy$
\[
\partial^{N}_{y}\tilde{\Psi}|_{\mathbb{R}}=(i\partial_{x})^{N}\tilde{\Psi}|_{\mathbb{R}}=(i\partial_{x})^{N}\Psi|_{\mathbb{R}}.
\]
In particular if $\langle x\rangle=(1+x^{2})^{1/2}$ then for any given  $N\geq
0$, a useful example of an almost analytic extension of $\Psi\in C^{\infty}_{0}(\mathbb{R})$ is given by
\[
\tilde{\Psi}(x+iy)=\Big(\sum_{m=0}^{N}\partial^{m}\Psi(x)(iy)^{m}/m!\Big)\tau(\frac{y}{\langle x \rangle}),
\]
where $\tau$ is a non-negative $C^{\infty}$ function such that
$\tau(s)=1$ if $|s|\leq 1$ and $\tau(s)=0$ if $|s|\geq 2$. For later
purposes, we also set
\[
\|\Psi\|_{N}\eqdef \sum_{m=0}^{N}\int_{\mathbb{R}}|\partial^{m}\Psi(x)|\langle x \rangle^{m-1}dx.
\]
\end{remark}
Our next lemma lets us deal with Lebesgue spaces.
\begin{lemma}\label{lemcontin}
Let $z\notin\mathbb{R}$ and $|\Im z|\lesssim |\Re z|$, then $\Delta_{D}$ satisfies 
\begin{equation}\label{hhh}
\|(z-\Delta_{D})^{-1}\|_{L^{p}(\Omega)\rightarrow
  L^{p}(\Omega)}\leq \frac{c}{|\Im z|}\left(\frac{|z|}{|\Im z|}\right)^{\alpha},\quad \forall z\notin\mathbb{R}
\end{equation}
for $1\leq  p\leq +\infty$, with a constant $c=c(p)>0$ 
and $\alpha=\alpha(n,p)>n|\frac 1 2 -\frac 1 p|$. 

Remark that, for all $h\in (0,1]$, the operator $h^{2}\Delta_{D}$
satisfies \eqref{hhh} with the same constants $c$ and $\alpha$ (this
is nothing but scale invariance).
\end{lemma}
For $p=2$ the proof of Lemma \ref{lemcontin} is trivial by multiplying the resolvent equation
$-\Delta_D u+z u=f$ by $\bar u$ and we get $\alpha=0$; however for $p\neq 2$ it requires a
non trivial argument which we postpone to the end of this Appendix.
\begin{corollary}
\label{cor}
For $N\geq \alpha+1$ the integral \eqref{integralpsi} is norm convergent and $\forall h\in (0,1]$
\begin{equation}\label{estimationpsi}
\|\Psi(-h^{2}\Delta_{D})\|_{L^{p}(\Omega)\rightarrow L^{p}(\Omega)}\leq c\|\Psi\|_{N+1},
\end{equation}
for some constant $c$ independent of $h$.
\end{corollary}
\begin{remark}
  Notice how the Mikhlin multiplier condition \eqref{eq:mm} on $\Psi$ does not
  imply boundedness of $\|\Psi\|_{N+1}$: we need extra decay at infinity.
\end{remark}
\emph{Proof:}
By scale invariance it is enough to prove \eqref{estimationpsi} for $h=1$. The integrand in \eqref{integralpsi} is norm continuous for $z\notin\mathbb{R}$. If we set
\[
U\eqdef \{z=x+iy|\langle x \rangle<|y|<2\langle x \rangle\},\quad V\eqdef \{z=x+iy| 0<|y|<2\langle x \rangle\},
\]
then the norm of the integrand is dominated by
\[
c\sum_{m=0}^{N}|\partial^{m}\Psi(x)|\frac{2^{m}}{m!}\langle x \rangle^{m-2}\|\partial\tau\|_{L^{\infty}([1,2])}1_{U}(x+iy)+
\]
\[
+c|\partial^{N+1}\Psi(x)|\frac{2^{N}}{N!}|y|^{N}\Big(\frac{\langle x \rangle}{|y|}\Big)^{\alpha}\|\tau\|_{L^{\infty}([0,2])}1_{V}(x+iy).
\]
Integrating with respect to $y$ for $N\geq \alpha+1$ yields the bound
\begin{multline*}
  \|\Psi(-\Delta_{D})\|_{L^{p}(\Omega)\rightarrow L^{p}(\Omega)}  \lesssim 
  \int_{\mathbb{R}}\Big(\sum_{m=0}^{N}|\partial^{m}\Psi(x)|\langle x \rangle^{m-1}+\\
 {}+|\partial^{N+1}\Psi(x)|\langle x \rangle^{N}\Big)dx=\|\Psi\|_{N+1}.
\end{multline*}
One may then prove that the operator $\Psi(-\Delta_{D})$, acting on
$L^p(\Omega)$, is independent of $N\geq 1+n/2$ and of the cut-off function
$\tau$ in the definition of $\tilde \Psi$, see \cite{Dav89}.

We now recall two lemma which will be useful when composing operators.
\begin{lemma}[Lemma 2.2.5,\cite{Dav89}]
\label{lemloc}
If $\Psi\in C^{\infty}_{0}(\mathbb{R})$ has support disjoint from the spectrum of $-h^{2}\Delta_{D}$ then $\Psi(-h^{2}\Delta_{D})=0$.
\end{lemma}
\begin{lemma}[Lemma 2.2.6, \cite{Dav89}]
\label{lemsupp}
If $\Psi_{1}$, $\Psi_{2}\in C^{\infty}_{0}(\mathbb{R})$, then $(\Psi_{1}\Psi_{2})(-h^{2}\Delta_{D})=\Psi_{1}(-h^{2}\Delta_{D})\Psi_{2}(-h^{2}\Delta_{D})$.
\end{lemma}


For the remaining part of the Appendix we prove the resolvent estimate \eqref{hhh} from Lemma
\ref{lemcontin}. If $\Re z>0$, this is nothing but a standard elliptic
estimate. The trouble comes with $\Re z <0$ and getting close to the
spectrum. In $\R^n$, one may evaluate directly the convolution
operator by proving its kernel to be in $L^1$: this follows from
$$
|z+|\xi|^2|^2 = \sin^2 \frac{(\pi-\theta)}{2} (|z|+|\xi|^2)^2+\cos^2\frac
{(\pi-\theta)}{2}(|\xi|^2-|z|)^2,\,\,\text{ with }\,\, z=|z|e^{i\theta},
$$
and a direct computation of $L^2$ norms of $\partial^\alpha
(z+|\xi|^2)^{-1}$. By reflection, one then extends this estimate to
the half-space case, with both Dirichlet and Neumann boundary
conditions. By localizing $L^p$ estimates close to the boundary and
flattening, one may then obtain the desired estimate \eqref{hhh}; such
an approach is carried out in \cite{3jap} in a greater generality
(systems of Laplace equations, mixed boundary conditions), at the
expense of fixing the angle $\theta$ and not tracking explicit
dependances on $|z|$ and $\theta$. While (relatively) elementary, such
a proof is, out of necessity, filled with lenghty calculations and
most certainly does not provide the sharpest constant. It is worth
noting, however, that it relies on standard elliptic techniques.

To keep in line with the parabolic approach, we present a short proof,
relying on the holomorphic nature of $S(w)$ in the half-plane $\Re
w>0$. Remark that by our $L^p$ bound on $S(t)$, $t\in \R_+$, the
trivial $L^2$ bound on $S(w)$, $\Re w\geq 0$, and Stein's parameter
version of complex interpolation, one may easily derive that $S(w)$ is
holomorphic in a sector around the positive real axis; but its angle
will narrow with large or small $p$. However the argument may be
refined and $S(w)$ was proved to be
holomorphic in the whole right half-plane in \cite{Ou}, using in a
crucial way the Gaussian nature of the heat kernel on domains
(\cite{DavHeat}).  This
was extented to more general settings in \cite{CCO}, where an explicit
bound is stated:
\begin{equation}
  \label{eq:CCO}
  \|S(w)\|_{L^p\rightarrow L^p} \leq C_\varepsilon \left(\frac{|w|}{|\Re
      w|}\right)^{n\left|\frac 1 2 -\frac 1 p\right|+\varepsilon}.
\end{equation}
Then \eqref{hhh} is a direct consequence of the following standard computation: recall the following formula, which is simply a
Laplace transform,
\begin{equation}
\label{laplace}
  (z -\Delta_D)^{-1}=\int_{L} e^{w\Delta_D-w z} dw,
\end{equation}
where $L$ can be chosen to be a half ray from the origin. Set
$z=r e^{i\theta}$, $w=\rho e^{i\phi}$, then
\begin{equation*}
  (z -\Delta_D)^{-1}=\int_0^{+\infty} e^{\rho \exp(i\phi)
    \Delta_D-r \rho \exp{i(\theta+\phi)}} d\rho.
\end{equation*}
Now, if $\Re z>0$, we may take $\phi=0$ and use
estimates for the semi-group $S(\rho)$. We would like to extend the range
to the $\Re z<0$ region, up to a thin sector around the negative real axis
($|\pi-\theta|<\epsilon$); getting close to the spectrum is required
if we want to define $\Psi(-\Delta_D)$ with
$\Psi \in C^\infty_0(]0,+\infty[)$. One picks $\phi$ such that
$2|\theta+\phi|<\pi$, which ensures a decaying exponential in
\ref{laplace}, provided we bound $S(w)$ in $L^p$. But the condition on
$\phi$ yields $|\phi|<\pi/2$, and the bound amounts to
the holomorphy of $S(w)$. The constant in \eqref{eq:CCO} translates
into a $(|z|/|\Im z|)^\alpha$ factor, while integration
over $\rho$ provides the remaining $1/|\Im z|$ in \eqref{hhh}. This concludes the
proof.
\bibliography{Besov}

\def\cprime{$'$}
\begin{thebibliography}{10}

\bibitem{3jap}
T.~Akiyama, H.~Kasai, Y.~Shibata, and M.~Tsutsumi.
\newblock On a resolvent estimate of a system of {L}aplace operators with
  perfect wall condition.
\newblock {\em Funkcial. Ekvac.}, 47(3):361--394, 2004.

\bibitem{Alex1}
Georgios~K. Alexopoulos.
\newblock ${L}^p$ bounds for spectral multipliers from {G}aussian estimates of
  the heat kernel.
\newblock unpublished manuscript, 1999.

\bibitem{Alex2}
Georgios~K. Alexopoulos.
\newblock Spectral multipliers for {M}arkov chains.
\newblock {\em J. Math. Soc. Japan}, 56(3):833--852, 2004.

\bibitem{Auscher}
Pascal Auscher.
\newblock On necessary and sufficient conditions for {$L\sp p$}-estimates of
  {R}iesz transforms associated to elliptic operators on {$\Bbb R\sp n$} and
  related estimates.
\newblock {\em Mem. Amer. Math. Soc.}, 186(871):xviii+75, 2007.

\bibitem{AMT}
Pascal Auscher, Alan McIntosh, and Philippe Tchamitchian.
\newblock Heat kernels of second order complex elliptic operators and
  applications.
\newblock {\em J. Funct. Anal.}, 152(1):22--73, 1998.

\bibitem{AT}
Pascal Auscher and Philippe Tchamitchian.
\newblock Square root problem for divergence operators and related topics.
\newblock {\em Ast\'erisque}, (249):viii+172, 1998.

\bibitem{BGT1}
N.~Burq, P.~G{\'e}rard, and N.~Tzvetkov.
\newblock Strichartz inequalities and the nonlinear {S}chr\"odinger equation on
  compact manifolds.
\newblock {\em Amer. J. Math.}, 126(3):569--605, 2004.

\bibitem{CCO}
Gilles Carron, Thierry Coulhon, and El-Maati Ouhabaz.
\newblock Gaussian estimates and {$L\sp p$}-boundedness of {R}iesz means.
\newblock {\em J. Evol. Equ.}, 2(3):299--317, 2002.

\bibitem{Daners}
Daniel Daners.
\newblock Heat kernel estimates for operators with boundary conditions.
\newblock {\em Math. Nachr.}, 217:13--41, 2000.

\bibitem{DavHeat}
E.~B. Davies.
\newblock {\em Heat kernels and spectral theory}, volume~92 of {\em Cambridge
  Tracts in Mathematics}.
\newblock Cambridge University Press, Cambridge, 1989.

\bibitem{DaviesJOT}
E.~B. Davies.
\newblock Pointwise bounds on the space and time derivatives of heat kernels.
\newblock {\em J. Operator Theory}, 21(2):367--378, 1989.

\bibitem{Dav89}
E.~B. Davies.
\newblock The functional calculus.
\newblock {\em J. London Math. Soc. (2)}, 52(1):166--176, 1995.

\bibitem{DaviesNGA}
E.~B. Davies.
\newblock Non-{G}aussian aspects of heat kernel behaviour.
\newblock {\em J. London Math. Soc. (2)}, 55(1):105--125, 1997.

\bibitem{Dynkin}
E.~M. Dyn{\cprime}kin.
\newblock An operator calculus based on the {C}auchy-{G}reen formula.
\newblock {\em Zap. Nau\v cn. Sem. Leningrad. Otdel. Mat. Inst. Steklov.
  (LOMI)}, 30:33--39, 1972.
\newblock Investigations on linear operators and the theory of functions, III.

\bibitem{HelSj}
B.~Helffer and J.~Sj{\"o}strand.
\newblock \'{E}quation de {S}chr\"odinger avec champ magn\'etique et \'equation
  de {H}arper.
\newblock In {\em Schr\"odinger operators ({S}\o nderborg, 1988)}, volume 345
  of {\em Lecture Notes in Phys.}, pages 118--197. Springer, Berlin, 1989.

\bibitem{OanaExtSchrod}
Oana Ivanovici.
\newblock On {S}chrodinger equation outside strictly convex obstacles, 2008.
\newblock {\tt arXiv:math.AP/0809.1060}.

\bibitem{GillesEvian}
Gilles Lebeau.
\newblock Estimation de dispersion pour les ondes dans un convexe.
\newblock In {\em Journ\'ees ``\'Equations aux D\'eriv\'ees Partielles''
  (Evian, 2006)}. 2006.

\bibitem{francis}
Francis Nier.
\newblock A variational formulation of {S}chr\"odinger-{P}oisson systems in
  dimension {$d\le 3$}.
\newblock {\em Comm. Partial Differential Equations}, 18(7-8):1125--1147, 1993.

\bibitem{Ou}
El-Maati Ouhabaz.
\newblock Gaussian estimates and holomorphy of semigroups.
\newblock {\em Proc. Amer. Math. Soc.}, 123(5):1465--1474, 1995.

\bibitem{PV}
Fabrice Planchon and Luis Vega.
\newblock Bilinear identities and applications, 2007.
\newblock to appear in Ann. Sci. E.N.S., {\tt arXiv:math.AP/0712.4076}.

\bibitem{SeeSog}
A.~Seeger and C.~D. Sogge.
\newblock On the boundedness of functions of (pseudo-) differential operators
  on compact manifolds.
\newblock {\em Duke Math. J.}, 59(3):709--736, 1989.

\bibitem{Stein}
Elias~M. Stein.
\newblock {\em Singular integrals and differentiability properties of
  functions}.
\newblock Princeton Mathematical Series, No. 30. Princeton University Press,
  Princeton, N.J., 1970.

\bibitem{DOS}
Xuan Thinh~Duong, El~Maati Ouhabaz, and Adam Sikora.
\newblock Plancherel-type estimates and sharp spectral multipliers.
\newblock {\em J. Funct. Anal.}, 196(2):443--485, 2002.

\bibitem{trev}
Fran{\c{c}}ois Tr{\`e}ves.
\newblock {\em Introduction to pseudodifferential and {F}ourier integral
  operators. {V}ol. 1}.
\newblock Plenum Press, New York, 1980.
\newblock Pseudodifferential operators, The University Series in Mathematics.

\bibitem{Triebel1}
H.~Triebel.
\newblock {\em Interpolation theory, function spaces, differential operators}.
\newblock VEB Deutscher Verlag der Wissenschaften, Berlin, 1978.

\bibitem{Triebel2}
Hans Triebel.
\newblock {\em Theory of function spaces}, volume~78 of {\em Monographs in
  Mathematics}.
\newblock Birkh\"auser Verlag, Basel, 1983.

\bibitem{Triebel3}
Hans Triebel.
\newblock {\em Theory of function spaces. {II}}, volume~84 of {\em Monographs
  in Mathematics}.
\newblock Birkh\"auser Verlag, Basel, 1992.

\end{thebibliography}
\end{document}